\documentclass[a4paper,12pt]{article}
\usepackage{latexsym}
\usepackage{amssymb}
\usepackage{amsfonts}
\usepackage{amsmath}

\newenvironment{@abssec}[1]{%
    \if@twocolumn

      \section*{#1}%
    \else

      \vspace{.05in}\footnotesize
      \parindent .2in
 {\upshape\bfseries #1. }\ignorespaces
    \fi}

    {\if@twocolumn\else\par\vspace{.1in}\fi}

\newenvironment{keywords}{\begin{@abssec}{\keywordsname}}{\end{@abssec}}

\newenvironment{AMS}{\begin{@abssec}{\AMSname}}{\end{@abssec}}

\newcommand\keywordsname{Key words}
\newcommand\AMSname{AMS subject classifications}
\newcommand\AMname{AMS subject classification}
\newtheorem{theorem}{Theorem}
 \newtheorem{lemma}[theorem]{Lemma}

\def\qed{\vbox{\hrule height0.6pt\hbox{%
  \vrule height1.3ex width0.6pt\hskip0.8ex
  \vrule width0.6pt}\hrule height0.6pt
 }}

\def\NN{\mathbb N}

\def\RR{\mathbb R}

\def\Om{\Omega}

\def\Ga{\Gamma}

\def\ka{\kappa}

\def\pa{\partial}

\title{Stationary isothermic surfaces \\ and some characterizations of the hyperplane\\ in the $N$-dimensional Euclidean space\thanks{This research was partially supported by a Grant-in-Aid
for
Scientific Research (B) ($\sharp$ 20340031) of
Japan Society for the Promotion of Science, and
by a
Grant of the Ital\-ian MURST.}}
\author{Rolando Magnanini\thanks{Dipartimento di Matematica U.~Dini,
Universit\` a di Firenze, viale Morgagni 67/A, 50134 Firenze, Italy.
({\tt magnanin@math.unifi.it}).}
  \quad and\quad Shigeru Sakaguchi\thanks{Department of Applied Mathematics,
Graduate School of Engineering, Hiroshima
University, Higashi-Hiroshima, 739-8527, Japan.
({\tt sakaguch@amath.hiroshima-u.ac.jp}).}}

\begin{document}

\maketitle

\begin{abstract}
We consider an entire graph $S$ of a continuous real function over $\mathbb R^{N-1}$ with $N\ge 3$. Let $\Omega$ be a domain in $\mathbb R^N$ with $S$ as a boundary. Consider in $\Omega$ the heat flow with initial temperature $0$ and boundary temperature $1.$ The problem we consider is to characterize $S$ in such a way that there exists a stationary isothermic surface in $\Omega$. We show that $S$ must be a hyperplane under some general conditions on $S$. This is related to Liouville or Bernstein-type theorems for some elliptic Monge-Amp\`ere-type equation.
\end{abstract}



\pagestyle{plain}
\thispagestyle{plain}
\markboth{R. MAGNANINI AND S. SAKAGUCHI}{Stationary isothermic surfaces}

\section{Introduction}

Let  $\Omega$ be a domain in $\mathbb R^N$ with $N \ge 3$, and let $u = u(x,t)$ be the unique bounded solution of the following problem for the heat equation: 
\begin{eqnarray}
&\partial_t u=\Delta u\ \ &\mbox{in }\ \Omega\times (0,+\infty),\label{heat}\\
&u=1\ \ &\mbox{on }\ \pa\Omega\times (0,+\infty),\label{dirichlet}\\
&u=0\ \ &\mbox{on }\ \Omega\times \{0\}.\label{initial}
\end{eqnarray}
The problem  we consider is to characterize the boundary $\partial\Omega$ in such a way that the solution $u$ has a stationary isothermic surface, say $\Gamma$. A hypersurface $\Gamma$ in $\Omega$ is said to be a {\it stationary isothermic surface} of $u$ if at each time $t $ the solution $u$ remains constant on $\Gamma$ ( a constant depending on $t$ ). 
It is easy to see that stationary isothermic surfaces occur when $\pa\Om$ and $\Ga$
are either parallel hyperplanes, concentric spheres, or coaxial spherical cylinders. 
The level surfaces of $u$ then are the so-called {\it isoparametric surfaces} whose
complete classification in Euclidean space was given by Levi-Civita \cite{LC} and Segre \cite{Seg}.

\par
Almost complete characterizations of the sphere have already been obtained by \cite{MS1, MS2} with the help of Aleksandrov's sphere theorem \cite{Alek}.
In \cite{MS2}, we also derived some characterizations of the hyperplane 
mainly based on geometrical arguments: under suitable global 
assumptions on $\partial\Omega,$ if $\Omega$ contains a stationary isothermic 
surface, then $\partial\Omega$ must be a hyperplane. In the present paper, we 
produce new results in this direction mainly based on partial differential 
equations techniques (in Section \ref{section3}, we compare them to the ones 
obtained in \cite{MS2}).
Assume that $\Omega$ satisfies the uniform exterior sphere condition and $\Omega$ is given by
\begin{equation}
\label{upper-half domain}
\Omega =\{\ x = (x^\prime,x_N) \in \mathbb R^N\ :\ x_N > \varphi(x^\prime)\  \},
\end{equation}
where $\varphi = \varphi(x^\prime)\ (x^\prime \in \mathbb R^{N-1})$ is a  continuous function on $\mathbb R^{N-1}$.
We recall that $\Omega$ satisfies the {\it uniform exterior sphere condition}
if there exists a number $r_0 >0$ such that  for every $\xi\in\pa\Omega$ 
there exists an open ball $B_{r_0}(y),$ centered at $y \in \mathbb R^N$ and with radius $r_0 > 0,$
satisfying $\overline{B_{r_0}(y)}\cap\overline{\Omega}=\{\xi\}.$ 
\par
We state our main result.

%%                      %%
%%  Theorem 1.1 begins. %%
%%                      %%

\begin{theorem}
\label{th:hyperplane}
Assume that there exists a stationary isothermic surface $\Gamma \subset \Omega$.
Then, under one of the following conditions {\rm (i), (ii),} and {\rm (iii)}, $\partial\Omega$ must be a hyperplane.
\begin{itemize}
\item[\rm(i)]\ $N =3;$
\item[\rm(ii)]\ $N \ge 4$ and $\varphi$ is globally Lipschitz continuous on $\mathbb R^{N-1};$
\item[\rm(iii)]\ $N \ge 4$ and there exists a non-empty open subset $A$ of $\partial\Omega$ such that on $A$ either  $H_{\partial\Omega}\ge 0$  or $\kappa_j \le 0\ $ for all $j = 1, \cdots, N-1.$
\end{itemize}
(Here $H_{\partial\Omega}$ and $\kappa_1, \cdots, \kappa_{N-1}$ denote the mean curvature of $\partial\Omega$ and the principal curvatures of $\partial\Omega$, respectively, with respect to the upward normal vector to $\partial\Omega.$)
\end{theorem}

\vskip 2ex
\noindent
{\bf Remark.} \ When $N=2$, this problem is easy. Since the curvature of  the curve $\partial\Omega$ is constant from (\ref{monge ampere}) in Lemma \ref{le:preliminary} in Section 2 of this paper, we see that $\partial\Omega$ must be a straight  line.
\par
Also, notice that, if $\varphi$ is either convex or concave, then (iii) is surely satisfied.

\setcounter{equation}{0}
\setcounter{theorem}{0}

\section{A proof of Theorem \ref{th:hyperplane} }
\label{section2}
The purpose of this section is to prove Theorem \ref{th:hyperplane}. 
Let $d = d(x)$ be the distance function defined by
\begin{equation}
\label{distance function}
d(x) = \mbox{ dist}(x, \partial\Omega),\quad x \in \Omega.
\end{equation}
We start with a lemma.
\begin{lemma}
\label{le:preliminary} The following assertions hold:
\begin{itemize}
\item[\rm(1)]\ $\Gamma = \{\ (x^\prime, \psi(x^\prime)) \in \mathbb R^N : x^\prime \in \mathbb R^{N-1}\ \}$ for some real analytic function $\psi = \psi(x^\prime)\ (x^\prime \in \mathbb R^{N-1})$;
\item[\rm(2)]\ There exists a number $R > 0$ such that $d(x) = R$ for every $x \in \Gamma$; 
\item[\rm(3)]\  $\varphi$ is real analytic, the mapping: 
$\partial\Omega \ni \xi \mapsto x(\xi) \equiv \xi + R\nu(\xi) \in \Gamma$ 
($\nu(\xi)$ denotes the upward unit normal vector to $\partial\Omega$ at $\xi \in \partial\Omega$) 
is a diffeomorphism, and $\partial\Omega$ and $\Gamma$ are parallel hypersurfaces at distance $R$;
\item[\rm(4)]\ the following inequality holds: for each $j =1, \cdots, N-1$
\begin{equation}
\label{bounds of curvatures}
-\frac 1{r_0} \le \kappa_j(\xi) < \frac 1R\  \mbox{ for every } \xi \in \partial\Omega,
\end{equation}
where $r_0 >0$ is the radius of the uniform exterior sphere for $\Omega$;
\item[\rm(5)]\ there exists a number $c > 0$ satisfying
\begin{equation}
\label{monge ampere}
\prod_{j=1}^{N-1} \left(\frac 1R-\kappa_j(\xi)\right) = c\quad\mbox{ for every } \xi \in \partial\Omega.
\end{equation}
\end{itemize}
\end{lemma}

\noindent
{\it Proof. }  The strong maximum principle implies that $\frac {\partial u}{\partial x_N} < 0$, and (1) holds. Since $\Gamma$ is stationary isothermic, (2) follows from  a result of Varadhan \cite{Va}:
$$
-\frac 1{\sqrt{s}} \log W(x,s) \to d(x)\ \mbox{ as } s \to \infty,
$$
where
\begin{equation}
\label{Laplace transform}
W(x,s) = s\int_0^\infty u(x,t) e^{-st} dt\ \mbox{  for } s >0.  
\end{equation}
The inequality $-\frac 1{r_0} \le \kappa_j(\xi)$ in (\ref{bounds of curvatures})  follows from the uniform exterior sphere condition for $\Omega$. See \cite[Lemma 2.2]{MS2} together with \cite[Lemma 3.1]{MS1} for the remaining claims. \qed

\vskip 3ex
With the help of Lemma \ref{le:preliminary}, we notice that $\varphi$ is an entire solution over $\mathbb R^{N-1}$ of the elliptic Monge-Amp\`ere-type equation (\ref{monge ampere}).
Thus, Theorem \ref{th:hyperplane} is related to Liouville or Bernstein-type theorems.

Let us proceed to the proof of Theorem \ref{th:hyperplane}.  Set
\begin{equation}
\label{middle surface}
\Gamma^* = \left\{\ x \in \Omega : d(x) = \frac R2\ \right\}.
\end{equation}
Denote by $\kappa^*_j$ and $\hat \kappa_j\ (j=1, \cdots, N-1)$ the principal curvatures of $\Gamma^*$ and $\Gamma$, respectively, with respect to the upward unit normal vectors. Then, the mean curvatures $H_{\Gamma^*}$ and $H_\Gamma$ of  $\Gamma^*$ and $\Gamma$ are given by
$$
H_{\Gamma^*} = \frac 1{N-1}\sum_{j=1}^{N-1} \kappa^*_j\ \mbox{ and }\ H_\Gamma =  \frac 1{N-1}\sum_{j=1}^{N-1} \hat \kappa_j,
$$
respectively. These principal curvatures have the following relationship: for each $j=1, \cdots, N-1,$
\begin{equation}
\label{relationship}
\kappa_j(\xi) = \frac {\kappa^*_j(\xi^*)}{1+\frac R2 \kappa^*_j(\xi^*)} = \frac {\hat \kappa_j(\hat\xi)}{1+R\hat \kappa_j(\hat \xi)}\ \mbox{ for any } \xi \in \partial\Omega,
\end{equation}
where $\xi^* = \xi + \frac R2\nu(\xi) \in \Gamma^*$ and $\hat \xi = \xi + R\nu(\xi) \in \Gamma$.
Let $\mu = c R^{N-1}.$ Then,  it follows from (\ref{monge ampere}) and (\ref{relationship}) that
\begin{equation}
\label{curvatures formulas}
\prod_{j=1}^{N-1}(1-R\kappa_j) = \mu,\  \prod_{j=1}^{N-1}(1+R\hat\kappa_j) = \frac 1\mu, \mbox{ and }\ \prod_{j=1}^{N-1}\frac {1-\frac R2\kappa^*_j}{1+\frac R2\kappa^*_j} =\mu.
\end{equation}
We distinguish three cases:
$$
{\rm ( I )}\ \mu > 1,\ {\rm ( II )}\ \mu < 1, \mbox{ and } {\rm ( III )}\ \mu = 1.
$$
Let us consider case (I) first. By the arithmetic-geometric mean inequality and the first equation in (\ref{curvatures formulas}) we have
$$
1-RH_{\partial\Omega} = \frac 1{N-1}\sum_{j=1}^{N-1} (1 -R\kappa_j) \ge \left\{ \prod_{j=1}^{N-1}(1-R\kappa_j)\right\}^{\frac 1{N-1}} = \mu^{\frac 1{N-1}} > 1.
$$
This shows that
\begin{equation}
\label{mean curvature inequality 1}
H_{\partial\Omega} \le -\frac 1R\left( \mu^{\frac 1{N-1}}-1\right) < 0.
\end{equation}
Since 
$$
(N-1)H_{\partial\Omega} = \mbox{ div}\left(\frac {\nabla\varphi}{\sqrt{1+|\nabla \varphi|^2}}\right)\ \mbox{ in } \mathbb R^{N-1}, 
$$
by using the divergence theorem we get a contradiction as in the proof of \cite[Theorem 3.3]{MS2}. 
In case (II), by the arithmetic-geometric mean inequality and the second equation in (\ref{curvatures formulas}) we have
$$
1+RH_{\Gamma} = \frac 1{N-1}\sum_{j=1}^{N-1} (1 +R\hat\kappa_j) \ge \left\{ \prod_{j=1}^{N-1}(1+R\hat\kappa_j)\right\}^{\frac 1{N-1}} = \mu^{-\frac 1{N-1}} > 1.
$$
This shows that
\begin{equation}
\label{mean curvature inequality 2}
H_{\Gamma} \ge \frac 1R\left( \mu^{-\frac 1{N-1}}-1\right) > 0,
\end{equation}
 which yields a contradiction similarly. 
 
 Thus, it remains to consider case (III). By  \eqref{mean curvature inequality 1}
and \eqref{mean curvature inequality 2}, we have
 \begin{equation}
 \label{super and sub}
 H_{\partial\Omega} \le 0 \le H_{\Gamma}.
 \end{equation}
 Let us consider case (i) of Theorem \ref{th:hyperplane} first. Since $N=3$ and $\mu=1$, it follows from the third equation of  (\ref{curvatures formulas}) that
 $$
2 H_{\Gamma^*} = \kappa^*_1+\kappa^*_2 = 0.
$$
We observe that $\Gamma^*$ is the entire graph of a function on $\mathbb R^2$. Therefore, by the Bernstein's theorem for the minimal surface equation,  $\Gamma^*$ must be a hyperplane. This gives the conclusion desired. (See \cite{GT, Giu} for Bernstein's theorem.)

Secondly, we consider case (iii) of Theorem \ref{th:hyperplane}.  Take any point $\xi \in A$.
If all the $\ka_j$'s are non-positive at $\xi,$ then they must vanish at $\xi,$
since $\prod\limits_{j=1}^{N-1}(1-R\kappa_j) = 1.$ 
On the other hand, we have that 
$$
1-RH_{\partial\Omega} = \frac 1{N-1}\sum_{j=1}^{N-1} (1 -R\kappa_j) \ge \left\{ \prod_{j=1}^{N-1}(1-R\kappa_j)\right\}^{\frac 1{N-1}} = 1;
$$
thus, if $H_{\pa\Om}\ge 0$ at  $\xi,$ then all the $\kappa_j$ must be equal
to each other and hence again they must vanish at $\xi$. Since $\xi \in A$ is arbitrary, we have
$$
\kappa_j \equiv 0\ \mbox{ on } A\ \mbox{ for every } j=1, \cdots, N-1,
$$
and hence $\varphi$ is affine on $A$.
Then by the analyticity of $\varphi$ we see that $\varphi$ is affine on the whole of $\mathbb R^{N-1}$.
This shows that $\partial\Omega$ must be a hyperplane.
\par
Thus it remains to consider case (ii) of Theorem \ref{th:hyperplane}. In this case, there exists a constant $L \ge 0$ satisfying
$$ 
\sup_{\mathbb R^{N-1}} |\nabla\varphi| =L < \infty.
$$
 Then, it follows from (1) and (3) of Lemma \ref{le:preliminary} that
\begin{equation}
\label{gradient bound}
\sup_{\mathbb R^{N-1}} |\nabla\psi|  = \sup_{\mathbb R^{N-1}} |\nabla\varphi| =L < \infty.
\end{equation}
Hence, in view of this and (3) of Lemma \ref{le:preliminary}, we can define a number $K^* > 0$ by
\begin{equation}
\label{critical number}
K^* = \inf\{ K > 0 : \psi\le \varphi+ K\ \mbox{ in } \mathbb R^{N-1} \}.
\end{equation}
Then we have
\begin{equation}
\label{useful inequality}
\varphi \le \psi \le h\ \mbox{ in } \mathbb R^{N-1},
\end{equation}
where $h:\mathbb R^{N-1}\to\RR$ is defined by
$$
h(x^\prime) = \varphi(x^\prime) + K^*\ \mbox{ for } x^\prime \in \mathbb R^{N-1}.
$$
\par
Moreover, by writing
$$
M(h) = \mbox{ div}\left(\frac {\nabla h}{\sqrt{1 +|\nabla h|^2}}\right) \mbox{ and }M(\psi) = \mbox{ div}\left(\frac {\nabla \psi}{\sqrt{1 +|\nabla \psi|^2}}\right), 
$$
from (\ref{super and sub}) and (\ref{useful inequality}) we have
\begin{equation}
\label{sub super method}
M(h) \le 0 \le M(\psi)\ \mbox{ and }\ \psi \le h\ \mbox{ in } \mathbb R^{N-1}.
\end{equation}
Hence,  the method of sub- and super-solutions with the help of (\ref{gradient bound}) yields that there exists $v \in C^\infty(\mathbb R^{N-1})$ satisfying
\begin{equation}
\label{middle solution}
M(v) = 0 \ \mbox{ and }\ \psi \le v \le h\ \mbox{ in } \mathbb R^{N-1}, \mbox{ and }\ \sup_{\mathbb R^{N-1}} |\nabla v| < \infty.
\end{equation}
\par
Indeed, take a sequence of balls $\{ B_n(0)\}_{n \in \mathbb N}$ in $\mathbb R^{N-1}$ and consider the boundary value problem for each $n \in \mathbb N$:
\begin{equation}
\label{Dirichlet problem}
M(v) = 0 \ \mbox{ in }\ B_n(0)\ \mbox{ and }\ v = \psi\ \mbox{ on } \partial B_n(0).
\end{equation}
By \cite[Theorem 16.9]{GT}, for each $n\in\NN$ there exists a $C^2$-function $v_n$ on $\overline{B_n(0)}$ solving problem (\ref{Dirichlet problem}). In view of (\ref{sub super method}), it then follows from the comparison principle that
\begin{equation}
\label{estimate form above and below}
\psi \le v_n \le  h\quad\mbox{ in }\ B_n(0)\quad \mbox{ for every } n \in \mathbb N.
\end{equation}
Therefore, with the help of the interior estimates for the minimal surface equation
(see \cite[Corollary 16.7]{GT}), we prove that $v_n$ belongs to $C^\infty(B_n(0))$ and for every  $\rho > 0$ and every $k \in \mathbb N$, the $C^k$ norms of $\{ v_n\}_{n > \rho}$ are bounded. 
In conclusion, the Cantor diagonal process together with Arzela-Ascoli theorem yields 
a solution $v \in C^\infty(\mathbb R^{N-1})$ of (\ref{middle solution}). 
It remains to show that $\nabla v$ is bounded in $\mathbb R^{N-1}$.
For this purpose, we  define a sequence of $C^\infty$ functions $\{ w_n \}$ on $\overline{B_1(0)}$ by
\begin{equation}
\label{scaling}
w_n(x^\prime) = \frac 1n\, v_n(nx^\prime)\ \mbox{ for }\ x^\prime \in B_1(0) \ \mbox{ and for every } n \in \mathbb N.
\end{equation}
Then, each $w_n$ satisfies
\begin{equation}
\label{Dirichlet problems on the unit ball}
M(w_n) = 0 \ \mbox{ in }\ B_1(0)\ \mbox{ and }\ w_n(x^\prime) = \frac 1n \psi(nx^\prime)\ \mbox{ for } x^\prime \in \partial B_1(0).
\end{equation}
\par
Since $ |(\nabla \psi)(nx^\prime)| \le L$, we have $ |\psi(nx^\prime)| \le  |\psi(0)| + n|x^\prime| L$. Therefore, it follows from the maximum principle that
$$
\max_{B_1(0)} |w_n| \le \max_{\partial B_1(0)} \frac 1n|\psi(nx^\prime)| \le  |\psi(0)| + L \ \mbox{ for every } n \in \mathbb N.
$$
Hence, by \cite[Corollary 16.7]{GT}, in particular there exists a constant $C$ satisfying
$$
|\nabla w_n(x^\prime)| \le C\ \mbox{ for every } x^\prime \in B_{\frac 12}(0) \ \mbox{ and for every } n \in \mathbb N.
$$
By observing that $\nabla w_n(x^\prime) = \left(\nabla v_n \right)(nx^\prime)$, we see that
$$
|\nabla v_n| \le C \ \mbox{ in } B_{n/2}(0)\ \mbox{ for every }n \in \mathbb N,
$$
 and hence
\begin{equation}
\label{bounded gradient}
|\nabla v | \le C\ \mbox{ in }\ \mathbb R^{N-1},
\end{equation}
which shows that the last claim in (\ref{middle solution}) holds.
Therefore, Moser's theorem \cite[Corollary, p. 591]{Mo} implies that $v$ is affine. We set $\eta = \nabla v \in \mathbb R^{N-1}$.
\par
On the other hand, by the definition of $K^*$ in (\ref{critical number}), there exists a sequence $\{ z_n \}$ in $\mathbb R^{N-1}$ satisfying
\begin{equation}
\label{critical sequence}
\lim_{n \to \infty} (h(z_n) -\psi(z_n)) = 0.
\end{equation}
Define a sequence of functions $\{ \varphi_n\}$ by
$$
\varphi_n(x^\prime) = h(x^\prime + z_n)-h(z_n) \ \left( = \varphi(x^\prime+z_n) - \varphi(z_n) \right).
$$
\par
Note that the principal curvatures $\kappa_1,\cdots, \kappa_{N-1}$ of $\partial\Omega$ are the eigenvalues of  the real symmetric matrix $G^{-\frac 12}BG^{-\frac 12},$ where the matrices $G$ and $B$ have entries
\begin{equation}
\label{first and second fundamental forms}
G_{ij} = \delta_{ij} + \frac {\partial\varphi}{\partial x_i}
\frac {\partial\varphi}{\partial x_j}
\ \mbox{ and }\ B_{ij} =\frac 1{\sqrt{1+|\nabla \varphi|^2}}
\frac {\partial^2\varphi}{\partial x_i\partial x_j}, 
\end{equation}
for $i,j=1,\dots, N-1,$ and $\delta_{ij}$ is Kronecker's symbol (see \cite[Proposition 3.1]{R}). 
\par
Then from (\ref{bounds of curvatures}) and (\ref{gradient bound}) we see that all the second derivatives of $\varphi$ are bounded in $\mathbb R^{N-1}$.
Hence we can conclude that there exists a subsequence $\{ \varphi_{n^\prime}\}$ of $\{ \varphi_n\}$ and a function $\varphi_\infty \in C^1(\mathbb R^{N-1})$ such that $\varphi_{n^\prime} \to \varphi_\infty$ in $C^1(\mathbb R^{N-1})$ as $n^\prime \to \infty$. Since $M(\varphi_n) \le 0$ in $\mathbb R^{N-1}$, we have that $M(\varphi_\infty) \le 0$ in $\mathbb R^{N-1}$ in the weak sense. Also, since $0 \le h(x^\prime + z_{n^\prime}) - v(x^\prime + z_{n^\prime})$ in $\mathbb R^{N-1}$, with the help of (\ref{critical sequence}),  letting $n^\prime \to \infty$ yields that 
$$
0 \le \varphi_\infty(x^\prime) - \eta\cdot x^\prime\ \mbox{ in } \mathbb R^{N-1}.
$$
\par
Consequently, we have
\begin{eqnarray*}
&&M(\varphi_\infty) \le 0  = M(\eta\cdot x^\prime)\ \mbox{ and } \varphi_\infty(x^\prime) \ge \eta\cdot x^\prime\ \mbox{ in } \mathbb R^{N-1}, 
\\ 
&&\mbox{ and } \varphi_\infty(0) =  0 = \eta \cdot 0.
\end{eqnarray*}
Hence, the strong comparison principle implies that $\varphi_\infty(x^\prime) \equiv \eta\cdot x^\prime$ in $\mathbb R^{N-1}$. Here we have used Theorem 10.7 together with Theorem 8.19 in \cite{GT}. Therefore we conclude that as $n \to \infty$
\begin{equation}
\label{c1 convergence1}
\varphi(x^\prime + z_n) - (v(x^\prime +z_n)-K^*) \to 0\ \mbox{ in } C^1(\mathbb R^{N-1}).
\end{equation}
\par
Similarly, we can obtain that as $n \to \infty$
\begin{equation}
\label{c1 convergence2}
v(x^\prime + z_n) - \psi(x^\prime+z_n) \to 0 \ \mbox{  in } C^1(\mathbb R^{N-1}).
\end{equation}
Indeed, it follows from (\ref{bounds of curvatures}) and (\ref{monge ampere}) that there exists a positive constant $\tau > 0$ such that for each $j =1, \cdots, N-1$
$$
-\frac 1{r_0} \le \kappa_j(\xi) \le \frac 1R-\tau\  \mbox{ for every } \xi \in \partial\Omega.
$$
Combining this with (\ref{relationship}) yields that all the principal curvatures $\hat \kappa_1, \cdots, \hat\kappa_{N-1}$ of $\Gamma$ are bounded.
Then,  in view of this fact and  the relationship between the function $\psi$ and  the principal curvatures $\hat \kappa_1, \cdots, \hat\kappa_{N-1}$,  from (\ref{gradient bound}) we see that all the second derivatives of $\psi$ are bounded in $\mathbb R^{N-1}$. Thus we can obtain (\ref{c1 convergence2}) by the same argument as in proving
(\ref{c1 convergence1}).

Therefore, it follows from (3) of Lemma \ref{le:preliminary}, (\ref{c1 convergence1}), and (\ref{c1 convergence2}) that the distance between two hyperplanes determined by two affine functions $v$ and $v-K^*$ must be $R$. Hence, since $v-K^* \le \varphi \le \psi \le v$ in $\mathbb R^{N-1}$,  we conclude that 
$$
\psi \equiv v\ \mbox{ and }\ \varphi \equiv v-K^* \ \mbox{ in } \mathbb R^{N-1},
$$
which shows that $\partial\Omega$ is a hyperplane.  \qed

\setcounter{equation}{0}
\setcounter{theorem}{0}

\section{Concluding remarks }
\label{section3}

Let us explain the relationship between Theorem \ref{th:hyperplane} and Theorems 3.2, 3.3, and 3.4 in \cite{MS2}. When $\mu = 1$, we have
$$
1+RH_{\Gamma} = \frac 1{N-1}\sum_{j=1}^{N-1} (1 +R\hat\kappa_j) \ge \left\{ \prod_{j=1}^{N-1}(1+R\hat\kappa_j)\right\}^{\frac 1{N-1}} = 1.
$$ 
Therefore, the assumption, $H_\Gamma \le 0$, of \cite[Theorem 3.2]{MS2}  implies that 
$\hat\kappa_j \equiv 0$ for every $j = 1, \cdots, N-1.$ This shows that $\Gamma$ is a hyperplane, and hence $\partial\Omega$ must be a hyperplane.
Thus, \cite[Theorem 3.2]{MS2} is contained in Theorem \ref{th:hyperplane} with its proof. In the case where $\Omega$ is given by (\ref{upper-half domain}), \cite[Theorem 3.3]{MS2} is contained in Theorem \ref{th:hyperplane} with condition (iii). Since \cite[Theorem 3.4]{MS2} does not assume the uniform exterior sphere condition for $\Omega$, it is independent of Theorem \ref{th:hyperplane}.

\end{document}